\documentclass[11pt]{amsart}
\usepackage{pstricks,pstcol,pst-plot}

\definecolor{Black}{cmyk}{0,0,0,1}
\definecolor{OrangeRed}{cmyk}{0,0.6,1,0}            % half magenta only, full yellow
\definecolor{DarkBlue}{cmyk}{1,1,0,0.20}
\definecolor{myblue}{rgb}{0.66,0.78,1.00}
\definecolor{Violet}{cmyk}{0.79,0.88,0,0}
\definecolor{Lavender}{cmyk}{0,0.48,0,0}

\parskip=\smallskipamount

\newtheorem{theorem}{Theorem}[section]
\newtheorem{lemma}[theorem]{Lemma}
\newtheorem{corollary}[theorem]{Corollary}

\theoremstyle{definition}

\newtheorem{example}[theorem]{Example}

\newtheorem{problem}[theorem]{Problem}
\newtheorem{remark}[theorem]{Remark}

\newcommand{\C}{\mathbb{C}}
\newcommand{\N}{\mathbb{N}}

\newcommand\hra{\hookrightarrow}

\def\di{\partial}
\def\bs{\backslash}
\def\e{\epsilon}

\numberwithin{equation}{section}

\begin{document}
\title[An interpolation theorem]
{An interpolation theorem \\ for proper holomorphic embeddings}
\author[F.\ Forstneri\v c, B.\ Ivarsson, F.\ Kutzschebauch and J.\ Prezelj]
{Franc Forstneri\v c, Bj\"orn Ivarsson, Frank Kutzschebauch \\ and Jasna Prezelj}
\address{F.\ Forstneri\v c and J.\ Prezelj: Institute of Mathematics, Physics and Mechanics, 
University of Ljubljana, Jadranska 19, 1000 Ljubljana, Slovenia}
\email{franc.forstneric@fmf.uni-lj.si \\ jasna.prezelj@fmf.uni-lj.si}
\address{B.\ Ivarsson and F.\ Kutzschebauch: Institute of Mathematics, University of Bern, Sidlerstrasse 5, CH-3012 Bern, Switzerland}
\email{bjoern.ivarsson@math.unibe.ch}
\email{Frank.Kutzschebauch@math.unibe.ch}  
\thanks{Forstneri\v c  and Prezelj supported by grants P1-0291 and J1-6173, Republic of Slovenia.}
\thanks{Kutzschebauch supported by Schweizerische Nationalfonds grant 
200021-107477/1.}
\thanks{Ivarsson supported by The Wenner-Gren Foundations.}

%
%    General info
%
\date{November 2, 2005} 

\begin{abstract}
Given a Stein manifold $X$ of dimension $n>1$, a discrete sequence $\{a_j\}\subset X$,
and a discrete sequence $\{b_j\}\subset \C^m$ where $m\ge N=\left[\frac{3n}{2}\right] + 1$,
there exists a proper holomorphic embedding $f\colon X\hra\C^m$ satisfying $f(a_j)=b_j$
for every $j=1,2,\ldots$.
\end{abstract}
\maketitle

\section{Introduction} 
It is known that every Stein manifold of dimension $n>1$ admits 
a proper holomorphic embedding in $\C^N$ with $N=\left[\frac{3n}{2}\right] + 1$,
and this $N$ is the smallest possible by the examples of Forster \cite{Fs1}.
The corresponding embedding theorem with $N$ replaced by 
$N'=\left[\frac{3n+1}{2}\right] + 1$ was proved  by 
Eliashberg and Gromov in \cite{EG} following an earlier
announcement in \cite{GE}. For even values of $n\in\N$ we have 
$N=N'$ and hence their result is the best possible; 
for odd values of $n$ the optimal result was obtained by Sch\"urmann \cite{Sch}, 
also for Stein spaces with bounded 
embedding dimension. A key ingredient in the known proofs of these results 
is the homotopy principle for holomorphic sections 
of elliptic submersions over Stein manifolds \cite{Gro}, \cite{FP2}. 

In this paper we prove the following embedding theorem with interpolation
on discrete sequences; for Stein spaces see theorem \ref{Steinspaces}.

%
%
%  Interpolation theorem for proper holo. embeddings
%
%
\begin{theorem} 
\label{Main}
Let $X$ be a Stein manifold of dimension $n>1$,
and let $\{a_j\}_{j\in\N}\subset X$ and $\{b_j\}_{j\in\N} \subset \C^m$
be discrete sequences without repetitions. 
If $m\ge N=\left[\frac{3n}{2}\right] +1$ then there exists a 
proper holomorphic embedding $f\colon X\hra \C^m$ satisfying 
\begin{equation} 
	f(a_j) = b_j \qquad (j=1,2,\ldots). 
\label{interpol}
\end{equation}
% $f(a_j)=b_j$ for $j=1,2,\ldots$. 
\end{theorem}

This result is optimal in all dimensions $n>1$ 
in view of Forster's examples \cite{Fs1}. 
For even values of $n\in\N$ theorem \ref{Main} coincides with 
a result of J.\ Prezelj to the effect that the conclusion holds with $N$ replaced 
by $N'=\left[\frac{3n+1}{2}\right] + 1$ (Theorem 1.1 (a) in \cite{Pr1}).
Our methods also give a different proof of Prezelj's result
to the effect that, under the assumptions of theorem \ref{Main} 
and with $m\ge \left[\frac{3n+1}{2}\right]$, there exists a 
proper holomorphic immersion $f\colon X\to \C^m$ satisfying 
(\ref{interpol}); see Theorem 1.1 (b) in \cite{Pr1}.

Prezelj obtained her results by carefully elaborating the constructions 
of Eliashberg and Gromov \cite{EG} and Sch\"urmann \cite{Sch}. 
It is not clear whether the method from \cite{Pr1} could be improved 
so as to give the optimal result also for odd values of $n$. 
We prove theorem \ref{Main} by
combining the known embedding theorems with methods of 
the theory of holomorphic automorphisms of Euclidean spaces.

If we increase the target dimension to $N\ge 2\dim X+1$ then
it is possible to extend any proper holomorphic embedding 
$Y\hra \C^N$ from an arbitrary closed complex submanifold $Y\subset X$ 
(not only a discrete set !) to a proper holomorphic embedding $X\hra \C^N$
(Acquistapace, Broglia and Tognoli \cite{ABT}; 
their proof follows closely those of
Bishop \cite{Bis} and Narasimhan \cite{Nar}).

Before proceeding, we recall that a discrete sequence 
$\{a_j\}_{j\in\N}$ in $\C^N$ is said to be 
{\em tame} in the sense of Rosay and Rudin \cite{RR}
if there exists a holomorphic automorphism of $\C^N$ which maps
$a_j$ to the point $e_j=(j,0,\ldots,0)$ for $j=1,2,\ldots$.
Several criteria for tameness can be found in \cite{RR};
for example, a sequence contained in a proper affine complex 
subspace of $\C^N$ is tame. 

Theorem \ref{Main} follows directly from the following two results.
The first one is seen by an inspection of the proofs in 
\cite{EG} and \cite{Sch} (see \S 3 below).
The second one is the main new result of this paper;  it has been proposed 
in \cite{BF1}, and it improves the result of \cite{Kut}. 

All sequences are assumed to be without repetition.

%
%
%  Eliashberg-Gromov-Schurmann
%
%
\begin{theorem}
\label{EGS}
{\rm (Eliashberg-Gromov-Sch\"urmann)}
Given a Stein manifold $X$ of dimension $n>1$ 
and a discrete sequence $\{a_j\}_{j\in \N} \subset X$,
there exists a proper holomorphic embedding
$f \colon  X\hra \C^N$ with $N=\left[\frac{3n}{2}\right] +1$
such that the sequence $\{f(a_j)\}_{j\in \N}$ is tame in $\C^N$. 
There also exists a proper holomorphic immersion 
$f\colon X\to \C^{[(3n+1)/2]}$ with the same property. 
\end{theorem}

%
%
%  Main theorem
%
%
\begin{theorem}
\label{Main1}
Let $N>1$, let $X$ be a closed, proper complex subvariety of $\C^N$,
and let $\{a_j\}_{j\in\N} \subset X$ be a discrete sequence 
which is tame in $\C^N$.  For every discrete sequence 
$\{b_j\}_{j\in\N} \subset \C^N$ there exist a domain 
$\Omega\subset \C^N$ containing $X$ and  a biholomorphic 
map $\Phi\colon\Omega\to \C^N$ onto $\C^N$
such that $\Phi(a_j)=b_j$ for $j=1,2,\ldots$.
\end{theorem}

Thus $X\to \Phi(X) \subset \C^N$ is another embedding 
of $X$ into $\C^N$ which interpolates the given sequences. 
In addition one can prescribe finite order jets
of $\Phi(X)$ at all points of the sequence which belong
to the regular locus of the subvariety (\S 2).
Note that $\Omega$ in theorem \ref{Main1} is a 
{\em Fatou-Bieberbach domain}. The fact that $\Phi(X)$ can be made to contain 
a given discrete sequence $\{b_j\}\subset\C^N$,
but without matching points, had been proved 
(for complex lines $\C\hra \C^2$) in \cite{FGR}, and in 
general in \cite{F1999}. Not surprisingly, the 
interpolation is considerably more difficult to achieve.

Since any discrete sequence contained in a proper  
{\em algebraic} subvariety of $\C^N$ is tame  \cite{RR},
theorem \ref{Main1} applies to all discrete
sequences $\{a_j\}\subset X$, $\{b_j\}\subset \C^N$ 
when $X$ is contained  in a proper algebraic subvariety
of $\C^N$.

Example \ref{ex1} below shows that theorem \ref{Main1} 
fails in general for non-tame sequences $\{a_j\}$.
The following problem of embedding with interpolation for a given
Stein manifold whose  embedding dimension is lower than the general 
dimension $N$ from theorem \ref{EGS}  therefore remains open.

\begin{problem}
\label{IP}
Let $X$ be a Stein manifold (or a Stein space) which admits a proper holomorphic 
embedding into $\C^m$ for some $m\in \N$. 
Given discrete sequences $\{a_j\}_{j\in\N}\subset X$ 
and $\{b_j\}_{j\in\N} \subset\C^m$ without repetitions, does there exist a 
proper holomorphic embedding $f\colon X\hookrightarrow \C^m$ satisfying
the interpolation condition (\ref{interpol})~?
\end{problem}

Since any discrete sequence in $\C^N=\C^N\times\{0\}\subset \C^{N+1}$ 
is tame in $\C^{N+1}$ \cite{RR}, theorem \ref{Main1} implies the 
following improvement of proposition 2.7 from \cite{Kut}
(adding only one extra dimension instead of two).

%
%
%   Corollary
%
%
\begin{corollary}
Let $X$ be a Stein space which admits a proper holomorphic embedding
into $\C^N$. If $m\ge N+1$ then for any two discrete sequences 
$\{a_j\}_{j\in\N}\subset X$ and $\{b_j\}_{j\in\N}\subset \C^m$ without repetitions
there exists a proper holomorphic embedding 
$f\colon X\hra \C^m$ satisfying  (\ref{interpol}).
\end{corollary}

The case $\dim X=1$, i.e., when $X$ is an open Riemann surface, 
is absent from the statement and discussion of theorem \ref{Main}.
The standard method fails when trying to embed such $X$ into $\C^2$
(it embeds into $\C^3$, also with interpolation on discrete sets 
\cite{ABT}, \cite{Bis}, \cite{Nar}).
For results in this direction see the survey \cite{F2003} 
and the recent papers of E.\ Forn\ae ss Wold 
(\cite{FW1}, \cite{FW2}) who showed in particular that every 
finitely connected planar domain embeds in $\C^2$, 
thereby extending the result of Globevnik and Stens\o nes \cite{GS}.

\begin{problem} 
For which open Riemann surfaces $X$ is problem \ref{IP}
solvable with $m=2$~? Is it solvable for every finitely connected planar domain~?
\end{problem}

Only two examples come to mind:
$X$ an algebraic curve in $\C^2$ when the result follows by applying 
theorem \ref{Main1}, and $X$ the unit disc when 
the interpolation theorem is due to Globevnik \cite{Glob}.

%
%
%  Section 2
%
%
\section{Proof of theorem \ref{Main1}}
We shall use the theory of holomorphic automorphisms of $\C^N$.
The precise result which we shall need is the following.

\begin{theorem}
\label{quote1}
{\rm (Buzzard and Forstneri\v c \cite{BF2}, Theorems 1.1 and 1.2.)}
Assume that $N>1$, $\{a_j\}$ and $\{a'_j\}$ are tame sequences
in $\C^N$, $K\subset\C^N$ is a compact, polynomially convex 
set contained in $\C^N\bs\{a_j\}$,
and $g$ is a holomorphic automorphism of $\C^N$ such that
$g(K)\subset \C^N\bs \{a'_j\}$. Then for every $\e>0$ 
there exists a holomorphic automorphism $\phi$ of $\C^N$
satisfying $\phi(a_j)=a'_j$ $(j=1,2,\ldots)$,  
$\sup_{z\in K} |\phi(z)-g(z)| <\e$, and
$\sup_{w\in g(K)} |\phi^{-1}(w)-g^{-1}(w)| <\e$. In addition 
one can prescribe finite order jets of $\phi$ at the
points $\{a_j\}$, and one can choose $\phi$ to exactly match 
$g$ up to a prescribed finite order at finitely many points of $K$. 
\end{theorem}

The statement concerning the approximation of $g^{-1}$
on $g(K)$ is a consequence of the approximation of $g$ on a
slightly larger polynomially convex set containing $K$ in its
interior, provided that $\e>0$ is sufficiently small.

The proof of theorem \ref{quote1} in \cite{BF2} relies upon 
the developments in \cite{AL}, \cite{FRo} and especially \cite{F1999}.
We shall use the special case of theorem \ref{quote1} 
when $g$ is the identity map and the 
sequence $\{a'_j\}$ differs from $\{a_j\}$ only in finitely many terms.
(Any modification of a tame sequence on finitely many terms 
is again tame.) The following lemma will provide the key step.

%
%
%  PROPOSITION - THE INDUCTIVE STEP
%
%
\begin{lemma}
\label{inductivestep}
Let $\{a_j\}\subset X\subset \C^N$ and $\{b_j\}\subset \C^N$ 
satisfy the hypotheses of theorem \ref{Main1}. 
Let $B\subset B'\subset \C^N$ be closed balls and $L=X\cap B'$. 
Assume that all points of the $\{b_j\}$ sequence which belong 
to $B\cup L$ coincide with the corresponding
points of the $\{a_j\}$ sequence, and all  remaining points 
of the $\{a_j\}$ sequence are contained in $X\bs L$.
Given $\e>0$ and a compact set $K\subset X$,
there exist a ball $B''\subset\C^N$ containing $B'$ 
($B''$ may be chosen as large as desired), a compact 
polynomially convex set $M\subset X$ with $L\cup K\subset M$,
and a holomorphic automorphism $\theta$ of $\C^N$
satisfying the following properties:
\begin{itemize}
\item[(i)]    $|\theta(z)-z|<\e$ for all $z\in B\cup L$,
\item[(ii)]   if $a_j\in M$ for some index $j$ then $\theta(a_j)=b_j\in B''$,
%(in particular, if $a_j \in L$ for some $j$ then $\theta(a_j)=a_j$), 
\item[(iii)]  if $b_j\in B'\bs (B\cup L)$ for some $j$ then $a_j \in M$ and $\theta(a_j)=b_j$, 
\item[(iv)]   $\theta(M) \subset {\rm Int} B''$, and
\item[(v)]    if $a_j\in X\bs M$ for some $j$ then $\theta(a_j)\in\C^N\bs B''$.
\end{itemize}
\end{lemma}

\begin{remark} 
\label{addendum}
If $\theta$ satisfies the conclusion of lemma \ref{inductivestep}
then the set
$$
	L'=\{z\in X\colon \theta(z)\in B''\}
$$
contains $M$ (and hence $K\cup L$), and $L'\bs M$ does not contain 
any points of the $\{a_j\}$ sequence (since the $\theta$-image
of any point $a_j\in X\bs M$ lies outside of $B''$ according to (v)). 
\end{remark}

\begin{proof}
An automorphism $\theta$ of $\C^N$ with the required properties 
will be constructed in two steps, $\theta=\psi\circ\phi$.

Since $X\cap B\subset L\subset X$ and the sets $B$ and $L$ 
are polynomially convex, $B\cup L$ is also polynomially convex 
(see e.g.\ Lemma 6.5 in \cite{F1999}, p.\ 111). 

By applying a preliminary automorphism of $\C^N$ which is very close 
to the identity map on $B\cup L$ we may assume that $X$ does not 
contain any points of the $\{b_j\}$ sequence, except 
those which coincide with the corresponding points $a_j\in X$. 
The same procedure will be repeated whenever necessary during 
later stages of the construction without mentioning it again.

Choose a pair of compact, polynomially convex neighborhoods 
$D_0\subset D\subset\C^N$ of $B\cup L$, with $D_0\subset {\rm Int} D$,
such that $D$ does not contain any additional 
points of the $\{a_j\}$ or the $\{b_j\}$ sequence. 
Choose $\e_0>0$ so small that 
$$
	{\rm dist}(B\cup L,\C^N \bs D_0)>\e_0, \quad 
        {\rm dist}(D_0,\C^N \bs D)>\e_0.
$$
By decreasing $\e>0$ if necessary we may assume $0<\e<\e_0$. 

Choose a compact polynomially convex set $M\subset X$
containing $K\cup (X\cap D)$ (and hence the set $L$),
and also containing all those points of the $\{a_j\}$
sequence for which the corresponding point $b_j$ is contained
in the ball $B'$. (Of course $M$ may also contain some additional
points of the $\{a_j\}$ sequence for which $b_j\in \C^N\bs B'$.)
Theorem \ref{quote1} furnishes an automorphism $\phi$ of $\C^N$ 
satisfying the following:
\begin{itemize}
\item[(a)]   $\sup_{z\in D} |\phi(z)-z|< \frac{\e}{2}$ and 
$\sup_{z\in D} |\phi^{-1}(z)-z|< \frac{\e}{2}$, 
\item[(b)]   $\phi(a_j)=b_j$ for all  $a_j\in M$, and
\item[(c)]   $\phi(a_j)=a_j$ for all $a_j\in X\bs M$. 
\end{itemize}

Condition (a) and the choice of $\e$ imply 
$\phi(D_0)\subset D$ and $\phi(\C^N\bs D)\cap D_0=\emptyset$,
and the latter condition also implies $\phi(X)\cap D_0 \subset \phi(M)$.
Since the sets $\phi(M)$ and $D_0$ are polynomially convex,
their union $\phi(M)\cup D_0$ is also polynomially convex 
(Lemma 6.5 in \cite{F1999}). 

Choose a large ball $B''\subset \C^N$ containing $\phi(M) \cup B'$. 
Theorem \ref{quote1} furnishes an automorphism 
$\psi$ of $\C^N$ satisfying the following:
\begin{itemize}
\item[(a')] $|\psi(z)-z|< \frac{\e}{2}$ when $z\in \phi(M)\cup D_0$,
\item[(b')] $\psi(\phi(a_j)) = \phi(a_j) =b_j$ for all $a_j\in M$, and
\item[(c')] $\psi(a_j)\in \C^N\bs B''$ for all $a_j\in X \bs M$.
\end{itemize}

We may also require that $\psi$ fixes all points $\phi (a_j)\in \phi(X)\bs B''$.
It is immediate that $\theta=\psi\circ\phi$ satisfies 
the conclusion of lemma \ref{inductivestep}.
\end{proof}

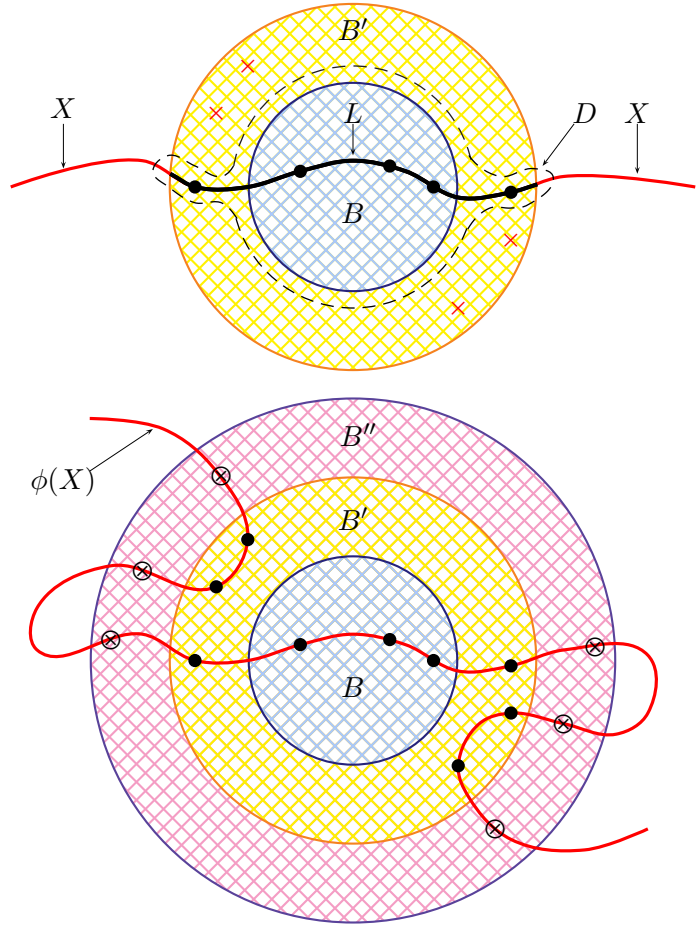
\begin{figure}[ht]
\psset{unit=0.7cm}

\begin{pspicture}(-7,-14)(7,4)
%
%
%  TOP PICTURE
%
%

\pscircle[linecolor=OrangeRed,fillstyle=crosshatch,hatchcolor=yellow](0,0){3.5}  %   the ball $B'$  
\rput(0,3){$B'$}						                 %    notation $B'$

\pscircle[linecolor=DarkBlue,fillstyle=crosshatch,hatchcolor=myblue](0,0){2}    %   the small ball $B$  
\rput(0,-0.5){$B$}					         	        %   notation $B$

\pscurve[linecolor=red,linewidth=1.2pt](-6.5,0)(-4,0.5)(-3,0)(-2,0)(-1,0.3)(0,0.5)
(1,0.3)(2,-0.2)(3,-0.1)(4,0.2)(6.5,0)                                           %   the analytic subvariety $X$

\psecurve[linecolor=black,linewidth=1.5pt](-4,0.5)(-3.47,0.25)(-3,0)(-2,0)(-1,0.3)(0,0.5)
(1,0.3)(2,-0.2)(3,-0.1)(3.5,0.05)(4,0.2)                                           %   the subset $L\subset X$

\psdots[dotsize=5pt](-3,0)(-1,0.3)(0.7,0.4)(1.53,0)(3,-0.1)

\psdots[dotstyle=x,dotscale=1.5,linecolor=red](-2.6,1.4)(-2,2.3)(3,-1)(2,-2.3)
%\psdots[dotstyle=+,dotangle=45,dotscale=1.5,linecolor=red](-2.6,1.4)(-2,2.3)(3,-1)(2,-2.3)
%\psdots[dotstyle=square,dotscale=1.2,linecolor=red](-2,2)

%
% The neighborhood $U$ of the set $B\cup L$, dashed.
%
\psarc[linestyle=dashed,linewidth=0.5pt,linecolor=Black](0,0){2.3}{13}{170}                     % upper dashed semicircle
\psecurve[linestyle=dashed,linewidth=0.5pt](2.1,1)(2.2,0.6)(2.8,0.2)(3.4,0.35)(3.8,0.2)
(3.3,-0.35)(3,-0.4)(2.6,-0.4)(2.18,-0.8)(2,-1.5)   						% right dashed curved
\psecurve[linestyle=dashed,linewidth=0.5pt](-2,2)(-2.27,0.4)(-2.5,0.2)(-3,0.4)(-3.5,0.6)
(-3.8,0.3)(-3.5,0)(-3,-0.3)(-2.5,-0.3)(-2.2,-0.6)(-2,-1)                                        % left dashed curve 
\psarc[linestyle=dashed,linewidth=0.5pt](0,0){2.3}{195}{340}                                    % bottom dashed semicircle

%
%   Notations $X$ and $D$, with accompanying arrows
%
\rput(-5.5,1.5){$X$}
\psline[linewidth=0.2pt]{->}(-5.5,1.2)(-5.5,0.35)
\rput(4.4,1.4){$D$}
\psline[linewidth=0.2pt]{->}(4.2,1.2)(3.6,0.4)
\rput(5.4,1.4){$X$}
\psline[linewidth=0.2pt]{->}(5.4,1.2)(5.4,0.25)

\rput(0,1.4){$L$}
\psline[linewidth=0.2pt]{->}(0,1.2)(0,0.6)

%
%
%  BOTTOM PICTURE
%
%

\pscircle[linecolor=Violet,fillstyle=crosshatch,hatchcolor=Lavender](0,-9){5}     %   the  largest ball $B''$  
\rput(0.1,-4.7){$B''$}						                 %    notation $B''$

\pscircle[linecolor=OrangeRed,fillstyle=crosshatch,hatchcolor=yellow](0,-9){3.5}  %   the ball $B'$  
\rput(0,-6.3){$B'$}						                 %    notation $B'$

\pscircle[linecolor=DarkBlue,fillstyle=crosshatch,hatchcolor=myblue](0,-9){2}    %   the small ball $B$  
\rput(0,-9.5){$B$}					         	          %   notation $B$

\pscurve[linecolor=red,linewidth=1.2pt] (-5,-4.4)(-3.6,-4.6)(-2.6,-5.4)(-2,-6.7)(-2.6,-7.6)(-4.6,-7.2)
(-6,-8)(-6,-8.8)(-4,-8.5)(-3,-9)(-2,-9)(-1,-8.7)(0,-8.5)(1,-8.7)(2,-9.2)(3,-9.1)(4,-8.8)
(5.6,-8.9)(5.6,-10)(5,-10.4)(4,-10.2)(3,-10)(2.4,-10.2)(2,-11)(2.2,-11.6)(3,-12.4)(4.4,-12.6)(5.6,-12.2)
%   the analytic subvariety $\phi(X)$

\rput(-5.5,-5.6){$\phi(X)$}
\psline[linewidth=0.2pt]{->}(-5,-5.4)(-3.8,-4.6)

\psecurve[linecolor=black,linewidth=1.5pt](-4,0.5)(-3.47,0.25)(-3,-0,3)(-2,0)(-1,0.3)(0,0.5)
(1,0.3)(2,-0.2)(3,-0.1)(3.5,0.05)(4,0.2)                                           %   the subset $L\subset X$

\psdots[dotsize=5pt](-3,-9)(-1,-8.7)(0.7,-8.6)(1.53,-9)(3,-9.1)(-2.6,-7.6)(-2,-6.7)(3,-10)(2,-11)
\psdots[dotstyle=otimes,dotscale=2](-2.5,-5.5)(-4,-7.3)(-4.6,-8.6)(4.6,-8.75)(4,-10.2)(2.7,-12.2)

\end{pspicture}
\caption{The proof of lemma \ref{inductivestep}.} 
\label{Fig1}
\end{figure}

The scheme of proof of lemma \ref{inductivestep}
is illustrated in fig.\ \ref{Fig1}. 
The first drawing shows the initial situation; the thick dots on $X$ 
indicate the points $b_j \in B\cup L$ which agree with the corresponding
points $a_j$, while the crosses indicate the remaining points $b_j\in B'$ 
which will be matched with the images of $a_j$ by applying 
the automorphism $\phi$. The second drawing shows 
the situation after the application of $\phi$: The large 
black dots in $\phi(X)\cap B'$ indicate the points $b_j=\phi(a_j) \in B'$, 
while the crossed dots on the subvariety $\phi(X)$ inside the set 
$B''\bs B'$ will be expelled from the ball $B''$ by the next 
automorphism $\psi$.

%
%
%  Figure 1:  Main Lemma 
%
%

\smallskip
{\em Proof of theorem \ref{Main1}.} 
Choose an exhaustion $K_1\subset K_2\subset \cdots \subset 
\bigcup_{j=1}^\infty K_j=X$ by compact sets. Fix a number $\e$
with $0<\e<1$. We shall inductively construct the following:
\begin{itemize}
\item[(a)] a sequence of holomorphic automorphisms $\Phi_k$ of $\C^N$ $(k\in\N)$,
\item[(b)] an exhaustion $L_1\subset L_2\subset \cdots\subset \bigcup_{j=1}^\infty L_j=X$
by compact, polynomially convex sets,
\item[(c)] a sequence of balls 
$B_1\subset B_2\subset\cdots \subset \bigcup_{j=1}^\infty B_j=\C^N$
centered at $0\in\C^N$ whose radii satisfy $r_{k+1}>r_k+1$
for $k=1,2,\ldots$,
\end{itemize}
such that the following hold for all $k=1,2,\ldots$
(conditions (iv) and (v) are vacuous for $k=1$):
\begin{itemize}
\item[(i)]  $\Phi_k(L_k)=\Phi_k(X) \cap B_{k+1}$,
\item[(ii)]  if $a_j\in L_k$ for some $j$ then $\Phi_k(a_j)=b_j$,
\item[(iii)] if $b_j \in \Phi_k(L_k)\cup B_k$ for some $j$ then $a_j \in L_k$ and $\Phi_k(a_j)=b_j$,
\item[(iv)]  $L_{k-1}\cup K_{k-1} \subset {\rm Int}L_{k}$, 
\item[(v)]  $|\Phi_{k}(z)-\Phi_{k-1}(z)|< \e\, 2^{-k}$ for all $z\in B_{k-1}\cup L_{k-1}$.
\end{itemize}

To begin we set $B_0=\emptyset$ and choose a pair of balls 
$B_1\subset B_2 \subset \C^N$ whose radii satisfy $r_2\ge r_1+1$. 
Theorem \ref{quote1} furnishes an automorphism $\Phi_1$ of $\C^N$ such that 
$\Phi_1(a_j)=b_j$ for all those (finitely many) indices $j$
for which $b_j\in B_2$, and $\Phi_1(a_j)\in\C^N\bs B_2$
for the remaining indices $j$. (Of course we only need to move
finitely many points of the $\{a_j\}$ sequence.)
Setting $L_1=\{z\in X\colon \Phi_1(z)\in B_2\}$,
the properties (i), (ii) and (iii) are satisfied for $k=1$
and the remaining two properties (iv), (v) are void.

Assume inductively that we have already found sets
$L_1,\ldots,L_k \subset X$, balls $B_1,\ldots,B_{k+1} \subset \C^N$
and automorphisms $\Phi_1,\ldots,\Phi_k$ such that 
(i)--(v) hold up to index $k$.
We now apply lemma \ref{inductivestep} with $B=B_k$,
$B'=B_{k+1}$, $X$ replaced by $X_k=\Phi_k(X) \subset\C^N$,
and $L=\Phi_k(L_k) \subset X_k$.
This gives us a compact polynomially convex set $M=M_k\subset X_k$
containing $\Phi_k(K_k\cup L_k)$, an automorphism $\theta=\theta_k$ of $\C^N$, 
and a ball $B''=B_{k+2} \subset \C^N$ of radius $r_{k+2} \ge r_{k+1}+1$ 
such that the conclusion of lemma \ref{inductivestep} holds.
In particular, $\theta_k(M_k) \subset B_{k+2}$, the 
interpolation condition is satisfied for all points 
$b_j \in \theta_k(M_k)\cup B_{k+1}$, and the remaining points
in the sequence $\{\Phi_k(a_j)\}_{j\in\N}$ are sent by 
$\theta_k$ out of the ball $B_{k+2}$. Setting 
$$ 
	\Phi_{k+1} =\theta_k\circ \Phi_k, \quad 
        L_{k+1}= \{z\in X\colon \Phi_{k+1}(z)\in B_{k+2}\}
$$ 
one easily checks that the properties (i)--(v) hold for the 
index $k+1$ as well. (Note that $L_{k+1}$ corresponds to the set $L'$
in remark \ref{addendum}). The induction may now continue.

Let $\Omega$ consist  of all points
$z\in\C^N$ for which the sequence $\{\Phi_k(z)\}_{k\in \N}$
remains bounded. Proposition 5.2 in \cite{F1999} (p.\ 108)
implies that $\lim_{k\to\infty} \Phi_k=\Phi$  exists
on $\Omega$, the convergence is uniform on compacts in
$\Omega$, and $\Phi\colon\Omega\to\C^N$ is a biholomorphic
map of $\Omega$ onto $\C^N$ (a Fatou-Bieberbach map).
In fact, $\Omega =\bigcup_{k=1}^\infty \Phi_k^{-1}(B_k)$
(Proposition 5.1 in \cite{F1999}). 
>From (v) we see that $X\subset \Omega$,
and properties (ii), (iii) imply that 
$\Phi(a_j)=b_j$ for all $j=1,2,\ldots$.
This completes the proof of theorem \ref{Main1}.

%
%  Example
%
\begin{example}
We show that theorem \ref{Main1} is not valid in general 
if $\{a_j\}$ is a non-tame sequence in $\C^N$.
Choose a sequence  $\{a_j\}_{j\in\N} \subset \C^N$ 
whose complement $\C^N \bs \{a_j\}_{j\in\N}$ is Eisenman $N$-hyperbolic 
\cite{Ka}, \cite{RR}. As already mentioned in the introduction, 
any complex subvariety $X\subset \C^N$ can be embedded in $\C^N$ so that 
its image contains a given sequence \cite{F1999}, and hence 
we may assume that $\{a_j\}_{j\in\N} \subset X$. 
Assume that theorem \ref{Main1} holds, i.e., there is 
a biholomorphic map $\Phi\colon \Omega\to \C^N$ from a domain
$\Omega\subset \C^N$ containing $X$ onto $\C^N$ satisfying $\Phi(a_j)=b_j$
for all $j=1,2,\ldots$. The set $\Omega \bs \{a_j\}_{j\in\N}$,
being contained in $\C^N \bs \{a_j\}_{j\in\N}$, is Eisenman $N$-hyperbolic,
and hence its $\Phi$-image $\C^N \bs \{b_j\}_{j\in\N}$ is Eisenman 
$N$-hyperbolic as well. But this is not true in general, for instance 
if the sequence $ \{b_j\}_{j\in\N}$ is tame in $\C^N$.
\label{ex1}
\end{example}

%
%
%  Section 3
%
%
\section{Embedding Stein spaces with interpolation}
We begin by indicating how theorem \ref{EGS} is obtained from 
Sch\"urmann's proof in \cite{Sch}.  

One begins by choosing a sufficiently generic 
almost proper holomorphic map $b\colon X\to \C^n$
with $n=\dim X$; this means that there are sequences
of compact special analytic polyhedra 
$K_1\subset K_2\subset\cdots\subset \bigcup_{j\in\N} K_j = X$ and
polydiscs $P_1\subset P_2\subset\cdots\subset \bigcup_{j\in\N} P_j=\C^n$
such that $b|_{K_j}\colon K_j\to P_j$ is a proper map 
sending the boundary $\di K_j$ to $\di P_j$ for every $j=1,2,\ldots$.
Such maps were first constructed by Bishop \cite{Bis} where
the reader can find more details; another source is 
Chapter VII in \cite{GR}. 

For a fixed $b$ as above one then constructs  a holomorphic map 
$g\colon X\to\C^{N-n}$ such that $f=(b,g)\colon X\hra \C^N$ 
is a proper holomorphic  embedding. The map $g$ is obtained
as the limit $g=\lim_{k\to\infty} g_k$ where the map $g_k\colon X\to\C^{N-n}$ 
accomplishes the job on $K_k$ and it approximates $g_{k-1}$ uniformly on $K_{k-1}$. 
The map $g$ has three tasks: 
to insure properness (this is done by choosing $|g_k|$ sufficiently 
large on $K_k\bs K_{k-1}$), to eliminate the kernel of the differential of $b$,
and to separate pairs of points which are not separated by $b$.
Such map can be found by the `elimination 
of singularities' method, due to Eliashberg and Gromov \cite{EG}, 
which proceeds by a finite induction over strata in a
suitable stratification of $X$. When extending the map from
one stratum to the next one uses the h-principle for sections
of elliptic submersions \cite{Gro}, \cite{FP2}. 
For the present purposes it is not necessary to understand this
method completely, and we refer the reader to \cite{EG} 
and \cite{Sch} for further details.

Suppose now that $\{a_j\}$ is a discrete sequence in $X$.
It is possible to choose the exhaustion of $X$ by special
analytic polyhedra $K_k$ as above such that $K_k\bs K_{k-1}$
contains at most one point of the sequence for each $k$.
Call this point $a_k$. When constructing the map
$g_k$ (which fulfills the relevant conditions on $K_k$)
it now suffices to require that the modulus of the last component
of the point $g_k(a_k)$ is sufficiently large;
it was already observed in \cite{Sch} and \cite{Pr1}
that this condition is easily built into the construction.
In this way we can achieve that the last components
of the sequence $\{g(a_j)\}_{j\in\N}$ form a discrete 
sequence (without repetitions) in $\C$. It follows from
standard methods (see e.g.\ \cite{RR}) that the sequence
$f(a_j)=(b(a_j),g(a_j)) \in\C^N$ is then tame.
This proves theorem \ref{EGS}.

Essentially the same proof applies if $X$ is a (reduced) Stein space with 
singularities and with bounded embedding dimension (Sch\"urmann  \cite{Sch}).
Let ${\rm Emb dim}_x X$ denote the local embedding dimension of $X$ at $x$,
that is, the smallest integer such that the germ of $X$ at $x$ 
embeds as a local closed complex subvariety of the Euclidean space 
of that dimension. Assume that 
$$
	q = {\rm Emb dim}\, X := \sup_{x\in X} {\rm Emb dim}_x X <+\infty. 
$$ 
Let $n(k)$ denote the dimension of the analytic set of points in $X$ 
at which $X$ has embedding dimension at least $k$. Set 
$$
	b\,'(X)=\max\{k+\left[n(k)/2\right] \colon k=0,\ldots,q\}.
$$ 
With this notation we have the following result,
extending theorem \ref{Main}.

\begin{theorem}
\label{Steinspaces}
Let $n>1$ and let $X$ be an $n$-dimensional Stein space 
of finite embedding dimension.  
Let $m\ge N =\max\{ \left[\frac{3n}{2}\right]+1,b\,'(X)\}$. 
Given discrete sequences $\{a_j\}\subset X$ and $\{b_j\} \subset \C^m$
without repetitions, there exists a proper holomorphic embedding 
$f\colon X\hra \C^m$ satisfying $f(a_j)=b_j$ for $j=1,2,\ldots$.
\end{theorem}

Theorem \ref{Steinspaces} is proved in the same way as theorem \ref{Main} by 
first embedding $X$ into $\C^m$ so that $\{a_j\}$
is mapped to a tame sequence in $\C^m$ (this is accomplished
by the modification of the proof in  \cite{Sch} described above), 
and subsequently applying theorem \ref{Main1}.

%
%
%  THANKS, THANKS AND THANKS
%
%
%
%\smallskip 
%\textit{Acknowledgements.} 

\bibliographystyle{amsplain}

\end{document}